\documentclass[12 pt]{amsart}

\usepackage{todonotes}
\usepackage{graphicx}
\usepackage{hyperref}

\newtheorem{thm}{Theorem}

\newtheorem{conj}[thm]{Conjecture}
\theoremstyle{definition}
\newtheorem{definition}[thm]{Definition}

\newtheorem{prop}[thm]{Proposition}
\newtheorem{example}[thm]{Example}
\theoremstyle{remark}
\newtheorem{remark}[thm]{Remark}

\newcommand\PP{{\mathbb{P}}}
\newcommand\Q{{\mathbb Q}}
\newcommand\N{{\mathbb{N}}}

\DeclareMathOperator{\inv}{inv}
\DeclareMathOperator{\asc}{asc}
\DeclareMathOperator{\rank}{rank}
\DeclareMathOperator{\DES}{DES}

\usepackage[top=1.2 in, bottom = 1.2 in, left = 1.2 in, right = 1.2 in]{geometry}
\begin{document}

\title[Chromatic quasisymmetric functions of directed graphs]
{ Chromatic quasisymmetric functions of directed graphs}
\author[Brittney Ellzey]{Brittney Ellzey$^1$}
\address{Department of Mathematics, University of Miami, Miami, FL}
\thanks{$^{1}$Supported in part by NSF Grant DMS 1202337}
\email{b.ellzey@math.miami.edu}
\date{November 14, 2016; revised April 1, 2017}

\begin{abstract}
{Chromatic quasisymmetric functions of labeled graphs were defined by Shareshian and Wachs as a refinement of Stanley's chromatic symmetric functions.  In this extended abstract, we consider an extension of their definition from labeled graphs to directed graphs, suggested by Richard Stanley.  We obtain an $F$-basis expansion of the chromatic quasisymmetric functions of all digraphs and a $p$-basis expansion for all \textit{symmetric} chromatic quasisymmetric functions of digraphs, extending work of Shareshian-Wachs and Athanasiadis.  We show that the chromatic quasisymmetric functions of proper circular arc digraphs are symmetric functions, which generalizes a result of Shareshian and Wachs on natural unit interval graphs.  The directed cycle on n vertices is contained in the class of proper circular arc digraphs, and we give a generating function for the $e$-basis expansion of the chromatic quasisymmetric function of the directed cycle, refining a result of Stanley for the undirected cycle.  We present a generalization of the Shareshian-Wachs refinement of the Stanley-Stembridge $e$-positivity conjecture.  }
\end{abstract}






\maketitle

\section{Introduction}

Let $G=(V,E)$ be a (simple) graph.  A proper coloring, $\kappa: V \rightarrow \PP$, of G is an assignment of positive integers, which we can think of as colors, to the vertices of G such that adjacent vertices have different colors; in other words, if $\{i,j\} \in E,$ then $\kappa(i) \neq \kappa(j).$  The chromatic polynomial of G, denoted $\chi_G(k),$ gives the number of proper colorings of G using $k$ colors.  Stanley \cite{CSF} defined a symmetric function refinement of the chromatic polynomial, called the chromatic symmetric function of a graph.  If we let the vertex set of G be $V=\{v_1, v_2,\cdots v_n\},$ then the chromatic symmetric function of G is defined as $$X_G({\bf x}) = \displaystyle \sum_{\kappa} x_{\kappa(v_1)}x_{\kappa(v_2)}\cdots x_{\kappa(v_n)}$$ where the sum ranges over all proper colorings, $\kappa$, of G and $\kappa(v_i)$ denotes the color of $v_i$.  The chromatic symmetric function of a graph refines the chromatic polynomial, because if we replace $x_1, x_2,\cdots,x_k$ with 1's and all other variables with 0's, then $X_G(1,1,\cdots,1,0,0,\cdots) = \chi_G(k).$

We can easily see that for any graph G, $X_G({\bf x}) \in \Lambda_\Q,$ where $\Lambda_\Q$ is the $\Q$-algebra of symmetric functions in the variables $x_1, x_2,\cdots$ with coefficients in $\Q.$  For any basis, $\{b_\lambda \mid \lambda \vdash n\},$ of $\Lambda_\Q$, we say that a symmetric function, $f \in \Lambda_\Q$ is $b$-positive if the expansion of the function in terms of the $b_\lambda$-basis has positive coefficients.  The symmetric function bases we focus on in this paper are the power sum symmetric function basis, $p_{\lambda},$ and the elementary symmetric function basis, $e_\lambda.$  We assume the reader is familiar with the basic theory of symmetric and quasisymmetric functions, which can be found in \cite{EC2}.

Stanley \cite{CSF} proves that $\omega X_G({\bf x})$ is $p$-positive for all graphs, $G$, where $\omega$ is the usual involution on $\Lambda_\Q.$ A long-standing conjecture on chromatic symmetric functions involves their $e$-positivity.  Recall that a poset is $(a+b)$-free if it has no induced poset that consists of a chain of length $a$ and a chain of length $b.$  The incomparability graph of a poset P, denoted $inc(P),$ is the graph whose vertices are the elements of P and whose edges correspond to pairs of incomparable elements of P.  The following conjecture is a generalization of a particular case of a conjecture of Stembridge on immanants \cite{Stem}.

\begin{conj} [Stanley-Stembridge \cite{CSF}]
Let P be a $(3+1)$-free poset.  Then $X_{inc(P)}({\bf x})$ is $e$-positive.
\end{conj}

For subsequent work on chromatic symmetric functions, see the references in \cite{CQSF}.


Shareshian and Wachs \cite{SW, CQSF} defined a quasisymmetric refinement of Stanley's chromatic symmetric function called the chromatic quasisymmetric function of a (labeled) graph.  Let $G = ([n],E)$ and let $\kappa:[n] \rightarrow \PP$ be a proper coloring of G.  Define the ascent number of $\kappa$ as $$\asc(\kappa) = |\{\{i,j\} \in E \mid i<j, \kappa(i)<\kappa(j)\}|.$$ The \textit{chromatic quasisymmetric function of G} is defined as $$X_G({\bf x},t) = \displaystyle \sum_{\kappa} t^{\asc(\kappa)}x_{\kappa(1)}x_{\kappa(2)} \cdots x_{\kappa(n)}$$ where $\kappa$ ranges over all proper colorings of G.  Notice that setting $t=1$ reduces this definition to Stanley's original chromatic symmetric function.

In the Shareshian-Wachs chromatic quasisymmetric function of a graph, it is not hard to see that the coefficient of $t^j$ for each $j \in \N$ is a quasisymmetric function; however, the coefficients do not have to be symmetric.  If $G$ is the path $1 - 2 - 3,$ then $X_G({\bf x},t)$ has symmetric coefficients, i.e. $X_G({\bf x},t) \in \Lambda_\Q[t]$, but if $G$ is the path $1 - 3 - 2,$ $X_G({\bf x},t)$ does not have symmetric coefficients (see \cite[Example 3.2]{CQSF}).  In general, $X_G({\bf x},t) \in QSym_\Q[t],$  where $QSym_\Q[t]$ is the ring of polynomials in t with coefficients in the ring of quasisymmetric functions in the variables $x_1, x_2, \cdots$ with coefficients in $\Q.$  


Shareshian and Wachs show that if G is natural unit interval graph (that is, a unit interval graph with a certain natural labeling), then $X_G({\bf x},t) \in \Lambda_\Q[t].$  For G a natural unit interval graph, they show that the coefficient of each power of t in $X_G({\bf x},t)$ is Schur-positive, and they conjecture that these coefficients are $e$-positive and $e$-unimodal.  In fact, Guay-Paquet \cite{GP} shows that if the Stanley-Stembridge conjecture holds for unit interval graphs, then the conjecture holds in general.  Hence the Shareshian-Wachs $e$-positivity conjecture implies the Stanley-Stembridge conjecture.  Shareshian and Wachs present a formula for the $e$-basis expansion of $X_{P_n}({\bf x},t),$ where $P_n$ is the path on n vertices with a natural labeling, showing that $X_{P_n}({\bf x},t)$ is $e$-positive. 

Shareshian and Wachs also conjectured a formula for the $p$-basis expansion of $\omega X_G({\bf x},t),$ where G is a natural unit interval order, which would imply that $\omega X_G({\bf x},t)$ is $p$-positive. Athanasiadis \cite{Ath} later proved this formula.

There is an important connection between chromatic quasisymmetric functions of natural unit interval graphs and Hessenberg varieties, which was conjectured by Shareshian and Wachs and was proven by Brosnan and Chow \cite{BrosChow} and later by Guay-Paquet \cite{GP2}. Clearman, Hyatt, Shelton, and Skandera \cite{Hecke} found an algebraic interpretation of chromatic quasisymmetric functions of natural unit interval graphs in terms of characters of type A Hecke algebras evaluated at Kazhdan-Lusztig basis elements. Recently, Haglund and Wilson \cite{HW} discovered a connection between chromatic quasisymmetric functions and Macdonald polynomials.

In this paper\footnote{i.e. in the full version of this paper \cite{me}}, we extend the work of Shareshian and Wachs by considering chromatic quasisymmetric functions of (simple) \textit {directed} graphs\footnote{See \hyperref[AP remark]{Remark \ref*{AP remark}}.}.  For notational convenience, we distinguish an undirected graph, G, from a digraph, $\overrightarrow{G},$ with an arrow.

\begin{definition} \label{def} Let $\overrightarrow{G}=(V,E)$ be a digraph.  Given a proper coloring, $\kappa:V \rightarrow \PP$ of $\overrightarrow{G}$, we define the ascent number of $\kappa$ as $$\asc(\kappa) = |\{(v_i, v_j) \in E \mid \kappa(v_i) < \kappa(v_j)\}|,$$ where $(v_i, v_j)$ is an edge directed from $v_i$ to $v_j.$  Then the \textit{chromatic quasisymmetric function of $\overrightarrow{G}$} is defined to be $$X_{\overrightarrow{G}}({\bf x},t) = \displaystyle \sum_{\kappa} t^{\asc(\kappa)} x_{\kappa(v_1)}x_{\kappa(v_2)} \cdots x_{\kappa(v_n)} $$ where the sum is over all proper colorings, $\kappa,$ of $\overrightarrow{G}$.   
\end{definition}

As with the Shareshian-Wachs chromatic quasisymmetric function, setting $t = 1$ gives Stanley's chromatic symmetric function.  We can easily see that for any digraph, $X_{\overrightarrow{G}}({\bf x},t) \in QSym_\Q[t].$  Notice that if we take a labeled graph $G = ([n],E)$ and make a digraph, $\overrightarrow{G},$ by orienting each edge from the vertex with the smaller label to the vertex with the larger label, then $X_G({\bf x},t) = X_{\overrightarrow{G}}({\bf x},t).$  In other words, this definition of the chromatic quasisymmetric function of a digraph is equivalent to the Shareshian-Wachs chromatic quasisymmetric function in the case of an acyclic digraph.

In this paper, we present an expansion of $\omega X_{\overrightarrow{G}}({\bf x},t)$ in Gessel's fundamental quasisymmetric basis with positive coefficients for every digraph, $\overrightarrow{G}.$  We determine a class of digraphs for which $X_{\overrightarrow{G}}({\bf x},t) \in \Lambda_\Q[t],$ namely proper circular arc digraphs.  A simple example of  a proper circular arc digraph is $\overrightarrow{C_n},$ the cycle on n vertices oriented cyclically.  If we turn natural unit interval graphs into directed unit interval graphs by orienting edges from smaller label to larger label, then these directed unit interval graphs are contained in the class of proper circular arc digraphs as well and hence our symmetry result generalizes the result of Shareshian and Wachs.

We present a $p$-positivity result for \textit{all} digraphs $\overrightarrow{G}$ such that $X_{\overrightarrow{G}}({\bf x},t) \in \Lambda_\Q[t],$ which does not reduce to the formula in the acyclic case conjectured by Shareshian-Wachs \cite{CQSF} and proved by Athanasiadis \cite{Ath}.  We give a factorization of the coefficients of $z_\lambda^{-1}p_{\lambda}$ in $\omega X_{\overrightarrow{C_n}}({\bf x},t).$  We also present a few results on $e$-positivity, including a generating function formula for $X_{\overrightarrow{C_n}}({\bf x},t),$ which is a $t$-analog of a result of Stanley \cite[Proposition 5.4]{CSF} and shows its $e$-positivity.  We present a generalization of the Shareshian-Wachs $e$-positivity conjecture for proper circular arc digraphs.  We also give a combinatorial interpretation of the coefficients in the elementary symmetric function expansion of the chromatic quasisymmetric functions of the cycle, oriented cyclically, and the path, oriented in one direction.


\begin{remark} \label{AP remark}
The idea of extending chromatic quasisymmetric functions to directed graphs was a suggestion made by Richard Stanley to the author after attending a talk on her work on the chromatic quasisymmetric function of the labeled cycle \cite{me}.

Subsequent to our work, Alexandersson and Panova \cite{AP} independently obtained the symmetry result of \hyperref[Sec 3]{Section \ref*{Sec 3}} and the results of \hyperref[Sec 4]{Section \ref*{Sec 4}}. However, their proof of \hyperref[cycle expansion]{Theorem~\ref*{cycle expansion}}, giving the $e$-expansion of $X_{\overrightarrow{C_n}}({\bf x},t)$, is very different from ours.

\end{remark}
\section{Expansion in the Fundamental Quasisymmetric Basis} \label{Sec 2}

For incomparability graphs of posets, Shareshian and Wachs give an expansion of $\omega X_G({\bf x},t)$ into Gessel's fundamental quasisymmetric basis, which shows that these $\omega X_{G}({\bf x},t)$ are $F$-positive.  We extend this result by presenting an $F$-basis expansion of $\omega X_{\overrightarrow{G}}({\bf x},t)$ for all digraphs, which shows that $\omega X_{\overrightarrow{G}}({\bf x},t)$ is $F$-positive for all digraphs.  In general our formula does not reduce to the formula of Shareshian and Wachs, so this gives another combinatorial description of the coefficients in the $F$-basis for incomparability graphs of posets.

Let $\overrightarrow{G}$ be a digraph and let $\sigma \in \mathfrak{S}_n.$  Define $$\inv_{\overrightarrow{G}}(\sigma) = |\{(i,j) \in E(\overrightarrow{G}) \mid \sigma^{-1}(j)< \sigma^{-1}(i)\}|,$$ i.e. the number of $(i,j)$ pairs such that j comes before i in $\sigma$ and there is a directed edge from i to j in $\overrightarrow{G}.$

Now let $G = ([n],E)$ be an undirected graph and let $\sigma = \sigma_1 \sigma_2 \cdots \sigma_n \in \mathfrak{S}_n$.  For each $x \in [n],$ define $\rank_{(G,\sigma)}(x)$ as the length of the longest subword $\sigma_{i_1}\sigma_{i_2}\cdots\sigma_{i_k}$ such that $i_1 < i_2 < \cdots < i_k$, $\sigma_{i_k} = x$, and for each $1 \leq j < k$, $\{\sigma_{i_j},\sigma_{i_{j+1}}\} \in E$.  We say $\sigma$ has a \textit{G-descent} at $i$ if either $\rank_{(G, \sigma)}(\sigma_i) > \rank_{(G, \sigma)}(\sigma_{i+1})$ or $\rank_{(G, \sigma)}(\sigma_i) = \rank_{(G, \sigma)}(\sigma_{i+1})$ and $\sigma_i > \sigma_{i+1}.$  Let $\DES_G(\sigma)$ be the set of G-descents of $\sigma.$

For example let $G = C_9,$ the cycle on 9 vertices labeled cyclically with $1, 2, ..., 9$ and let $\sigma = 234658971 \in \mathfrak{S}_9.$ For $x \in [9],$ if $\rank_{(G, \sigma)}(x) = 1,$ then $x = 2, 6, 8.$  If $\rank_{(G, \sigma)}(x) = 2,$ then $x = 3, 7, 9.$ If $\rank_{(G, \sigma)}(x) = 3,$ then $x = 1,4.$ If $\rank_{(G, \sigma)}(x) = 4,$ then $x = 5.$  From this we see that $\DES_G(\sigma) = \{3, 5, 7\}.$

\begin{thm} \label{F-basis thm}
Let $\overrightarrow{G} = ([n],E)$ be any directed graph.  Then $$\omega X_{\overrightarrow{G}}({\bf x},t) = \displaystyle \sum_{\sigma \in \mathfrak{S}_n} F_{n, \DES_G(\sigma)}({\bf x})t^{\inv_{\overrightarrow{G}}(\sigma)}$$ where $F_{n,S}({\bf x})$ is Gessel's fundamental quasisymmetric function and G is the underlying undirected graph of $\overrightarrow{G}.$
\end{thm}

Note that our formula requires that $\overrightarrow{G}$ be labeled with $[n]$; however, the expansion of $X_{\overrightarrow{G}}({\bf x},t)$ in the $F$-basis is the same for any choice of labeling.
\section{Expansion in the Power Sum Symmetric Function Basis} \label{p}
In \cite{CSF}, Stanley showed that for \textit{any} graph $G,$ $\omega X_G({\bf x})$ is $p$-positive.  In \cite{CQSF}, Shareshian and Wachs conjectured a formula for the $p$-expansion coefficients that established that $\omega X_G({\bf x},t)$ is $p$-positive for any natural unit interval graph, $G$, and in \cite{Ath}, Athanasiadis proved their conjecture. Here we present a $p$-positivity result for all digraphs whose chromatic quasisymmetric functions have symmetric coefficients.

We give a formula for the coefficients of each $p_{\lambda}$ in the $p$-basis expansion of $\omega X_{\overrightarrow{G}}({\bf x},t).$ We can assume without loss of generality that the vertex set of $\overrightarrow{G}$ is $[n].$  We want to define a set of permutations, $N_{G,\lambda},$ for every undirected graph $G = ([n],E)$ and every partition $\lambda$ of n.  


For any word $w = w_1 w_2 \cdots w_k$ with distinct letters in $[n],$ we say a letter $w_j$ with $j>1$ is {\textit{G-isolated}} if for all $i<j$, there is no edge between $w_i$ and $w_j$ in $G.$

We define $N_{G,\lambda}$ as follows.  If $\sigma \in \mathfrak{S}_n,$ break $\sigma$ up into contiguous segments of sizes $\lambda_1, \lambda_2, \cdots, \lambda_k$ (in order) and call these pieces $\alpha_1, \alpha_2, \cdots \alpha_k.$  Then $\sigma \in N_{G,\lambda}$ means that each segment $\alpha_i$ has no G-isolated letters and does not contain any of the $G$-descents of $\sigma$.  So every permutation is in $N_{G, (1)^n},$ and the permutations of $N_{G,(n)}$ are in bijection with acyclic orientations of G with a unique sink.

In the last section, we determined that for $G = C_9$ and $\sigma = 234658971,$ $\DES_G(\sigma) = \{3,5,7\}.$  If we let $\lambda = (3, 2, 2, 1, 1),$ then $\alpha_1 = 234, \alpha_2 = 65, \alpha_3 = 89, \alpha_4 = 7,$ and $\alpha_5 = 1,$ and $\sigma \in N_{G, \lambda};$ however if $\mu = (3, 2, 2, 2),$ then $1$ is an isolated vertex of $\alpha_4 = 71,$ so $\sigma \notin N_{G, \mu}.$ 

\begin{thm}
Let $\overrightarrow{G} = ([n],E)$ be a digraph such that $X_{\overrightarrow{G}}({\bf x},t) \in \Lambda_Q[t].$ Then $$\omega X_{\overrightarrow{G}}({\bf x},t) = \displaystyle \sum_{\lambda \vdash n} z_{\lambda}^{-1} p_{\lambda} \displaystyle \sum_{\sigma \in N_{G,\lambda}}t^{\inv_{\overrightarrow{G}}(\sigma)},$$ where $G$ is the underlying undirected graph of $\overrightarrow{G}.$  Consequently, $\omega X_{\overrightarrow{G}}({\bf x},t)$ is p-positive.
\end{thm}

We prove this by first expressing $\omega X_G({\bf x},t)$ in terms of the fundamental quasisymmetric basis, as shown in \hyperref[F-basis thm]{Theorem \ref*{F-basis thm}}, and then extending the technique used by Athanasiadis \cite{Ath}, which involves the Adin-Roichman \cite{Adin} formula for symmetric group representations on Schur-positive sets, to prove the acyclic version of this result.  We point out that our result does not reduce, in an obvious way, to the Athanasiadis-Shareshian-Wachs formula in the acyclic case.  It gives a new formula for this case.

In \cite[Proposition 7.8]{CQSF}, Shareshian and Wachs showed that when $\overrightarrow{G}$ is acyclic, the coefficient of each $z_\lambda^{-1}p_{\lambda}$ in $\omega X_{\overrightarrow{G}}({\bf x},t)$ factors nicely. Though the coefficients of $\omega X_{\overrightarrow{G}}({\bf x},t)$ do not generally factor in the cyclic case, the coefficient of each $z_\lambda^{-1}p_{\lambda}$ in $\omega X_{\overrightarrow{C_n}}({\bf x},t)$ does have a nice factorization involving the Eulerian polynomials.  

\begin{thm}
Let $\overrightarrow{C_n}$ be the cycle on n vertices directed cyclically and let $\lambda = (\lambda_1, \lambda_2, \cdots, \lambda_k)$ be a partition of n.  If $k \geq 2$, then $$\displaystyle \sum_{\sigma \in N_{C_n, \lambda}}t^{\inv_{\overrightarrow{C_n}}(\sigma)} = nt A_{k-1}(t) \displaystyle \prod_{i = 1}^k [\lambda_i]_t,$$ where $[n]_t = 1 + t + \cdots + t^{n-1}$ and $A_k(t)$ is the Eulerian polynomial.  In the case that $\lambda = (n),$ we have $$\displaystyle \sum_{\sigma \in N_{C_n, (n)}}t^{\inv_{\overrightarrow{C_n}}(\sigma)} = nt[n-1]_t,$$ Hence the coefficient of $\frac{1}{n}p_n$ in $\omega X_{\overrightarrow{C_n}}({\bf x},t)$ is $nt[n-1]_t$ and for all other $\lambda \vdash n$, the coefficient of $z_{\lambda}^{-1}p_\lambda$ in $\omega X_{\overrightarrow{C_n}}({\bf x},t)$ is $nt  A_{k-1}(t) \displaystyle \prod_{i = 1}^k [\lambda_i]_t.$ 
\end{thm}
\section{Graphs and Symmetry} \label{Sec 3}

In this section, we discuss a class of digraphs, $\overrightarrow{G}$, such that $X_{\overrightarrow{G}}({\bf x},t)$ is symmetric.  An oriented graph\footnote{An oriented graph is a digraph with no bidirected edges.} will be called a {\textit{$\{\overrightarrow{K_{12}}, \overrightarrow{K_{21}}\}$-free digraph} if it avoids the induced digraphs $\overrightarrow{K_{12}}$ and $\overrightarrow{K_{21}},$ shown below.

\begin{center}
\includegraphics[scale = 0.5]{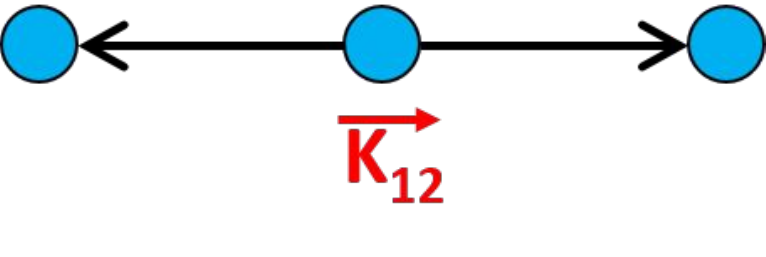} \hspace{0.5 in} \includegraphics[scale=0.5]{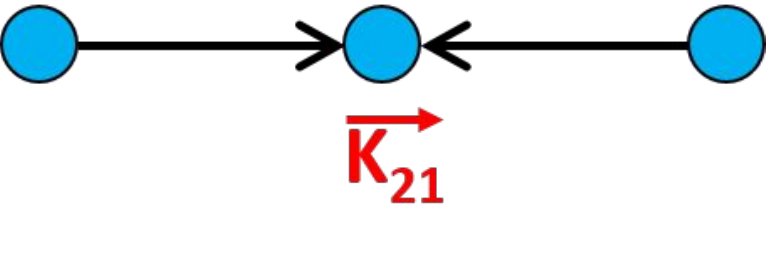}
\vspace{-0.2 in}
\end{center}

The most well-known class of graphs discussed in this paper are \textit{interval graphs}.  Given a collection of intervals on the real line, we can associate them with a graph by letting each interval correspond to a vertex and each edge correspond to a pair of overlapping intervals.  \textit{Proper interval graphs} are interval graphs in which no interval properly contains another.  \textit{Unit interval graphs} are interval graphs in which each interval has unit length. 

\begin{prop}[Roberts \cite{Roberts}, Skrien \cite{Skrien}]
Let G be a graph. Then the following statements are equivalent:

1. G is a proper interval graph.

2. G is a unit interval graph.

3. G admits an acyclic orientation that makes it a $\{\overrightarrow{K_{12}}, \overrightarrow{K_{21}}\}$-free digraph.

\end{prop}

The equivalence of (1) and (2) was shown by Roberts \cite{Roberts} and the equivalence of (1) and (3) was shown by Skrien \cite{Skrien}.  

In \cite{CQSF}, Shareshian and Wachs prove that chromatic quasisymmetric functions of labeled graphs have symmetric coefficients for natural unit interval graphs, which are unit interval graphs with a specific labeling. For example, the path $1 - 2 - 3$ is a natural unit interval graph, but the paths $1 - 3 - 2$ and $2 - 1 - 3$ are not.  Note that the digraph associated with $2 - 1 - 3$ is $\overrightarrow{K_{12}}$ and the digraph associated with $1 - 3 - 2$ is  $\overrightarrow{K_{21}}.$ For the remainder of the paper, we will use the term \textit{unit interval digraphs} to refer to natural unit interval graphs viewed as digraphs.  It turns out that this is exactly the class of acyclic $\{\overrightarrow{K_{12}}, \overrightarrow{K_{21}}\}$-free digraphs.  

Now let us look at the circular analog of these graph classes.  The circular analog of interval graphs is the class of \textit{circular arc graphs}.  If we have a collection of arcs on a circle, we can associate a graph to this collection by allowing each arc to correspond to a vertex and each edge to correspond to a pair of overlapping arcs.  \textit{Proper circular arc graphs} are circular arc graphs where no arc properly contains another.  
 
\begin{prop}[Skrien \cite{Skrien}]
Let G be a connected graph. Then the following statements are equivalent:

1. G is a proper circular arc graph.

2. G admits an orientation that makes it a $\{\overrightarrow{K_{12}}, \overrightarrow{K_{21}}\}$-free digraph.
\end{prop}

For the remainder of this paper, we will refer to $\{\overrightarrow{K_{12}}, \overrightarrow{K_{21}}\}$-free digraphs as \textit{proper circular arc digraphs}, because this proposition shows that the underlying undirected graph of each of the connected components of a $\{\overrightarrow{K_{12}}, \overrightarrow{K_{21}}\}$-free digraph is a proper circular arc graph.

The smallest digraphs that do not have symmetric coefficients are $\overrightarrow{K_{12}}$ and $\overrightarrow{K_{21}}.$  By the work of Shareshian and Wachs, we know that chromatic quasisymmetric functions of unit interval digraphs, or acyclic $\{\overrightarrow{K_{12}}, \overrightarrow{K_{21}}\}$-free digraphs, have symmetric coefficients.  We have the following generalization of this result.  The proof closely follows the proof of the Shareshian-Wachs symmetry result \cite[Theorem 4.5]{CQSF}.

\begin{thm} \label{sym} 
Let $\overrightarrow{G}$ be a proper circular arc digraph and let $X_{\overrightarrow{G}}({\bf x},t)$ be the chromatic quasisymmetric function associated with the digraph $\overrightarrow{G}$.  Then $X_{\overrightarrow{G}}({\bf x},t) \in \Lambda_{\Q}[t].$
\end{thm}

Note that the converse of this statement is not true.  In fact, the cycle with one edge directed backwards is in $\Lambda_\Q[t]$; however, for the rest of this paper, we will focus on the chromatic quasisymmetric functions of proper circular arc digraphs.

\section{Expansion in the Elementary Symmetric Function Basis}\label{Sec 4}

The Shareshian-Wachs conjecture stated in terms of digraphs says that the chromatic quasisymmetric functions of unit interval digraphs are $e$-positive and $e$-unimodal, where we call the palindromic\footnote{This is shown to be a palindromic polynomial in \cite{CQSF}.} polynomial $X_G({\bf x},t) = \sum_{j = 0}^m a_j({\bf x})t^j$ \textit{e-unimodal} if $a_{j+1}({\bf x}) - a_j({\bf x})$ is $e$-positive for $0 \leq j < \frac{m-1}{2}.$  We present a generalized version of this conjecture\footnote{An equivalent form of this conjecture is also noted in \cite{AP}.}.

\begin{conj}\label{conj}
Let $\overrightarrow{G}$ be a proper circular arc digraph.  Then the palindromic\footnote{This is shown to be a palindromic polynomial in \cite{me}.} polynomial $X_{\overrightarrow{G}}({\bf x},t)$ is $e$-positive and $e$-unimodal.
\end{conj}

For $r, n \in \N$ such that $1 \leq r \leq n$ define $G_{n,r} = ([n],E)$ to be the graph on $[n]$ with $\{i,j\} \in E$ if $0 < |i - j| < r.$  For example, $G_{n,1}$ is the graph on $[n]$ with no edges, and $G_{n,2}$ is the path on $[n]$, where consecutive labels are adjacent.  It is not difficult to see that $X_{G_{n,1}}({\bf x},t) = e_1^n.$ Shareshian and Wachs proved the conjecture for $G_{n, n-1}$ and $G_{n,n-2}$ \cite[Corollaries 8.2, 8.3, 8.4]{CQSF} as well as for the path, $G_{n,2}$ \cite[Theorem 7.2]{Eul}, and the complete graph, $G_{n,n}$ for all $n,$ and they tested it for all $G_{n,r}$ with $n \leq 8$ and $1 \leq r \leq n,$ which also implies that this conjecture holds for these $G_{n,r}$ graphs viewed as digraphs rather than labeled graphs.

For $r, n \in \N$ such that $1 \leq r \leq \left \lceil{\frac{n}{2}}\right \rceil,$ define the directed circular analog $\overrightarrow{G}_{n,r}^* = ([n],E)$ to be the digraph on $[n]$ with $(i,j) \in E$ if $0 <(j - i)  \pmod n <r.$ For example, $\overrightarrow{G}_{n,2}^*$ is the cycle on n vertices with edges directed cyclically.  We show in the next theorem that $X_{\overrightarrow{G}_{n,2}^*}({\bf x},t)$ is $e$-positive and $e$-unimodal.  We tested the $e$-positivity and $e$-unimodality of $\overrightarrow{G}_{n,r}^*$ for $n \leq 9$ and $r \leq \left \lceil{\frac{n}{2}}\right \rceil.$

Note that \hyperref[conj]{Conjecture \ref*{conj}} would imply the Shareshian-Wachs conjecture, which in turn would imply the Stanley-Stembridge conjecture.  In addition, Stanley \cite{CSF} mentions a class of graphs called circular indifference graphs, which turns out to be equivalent to the class of proper circular arc graphs, as a possible class of graphs with $e$-positive chromatic symmetric functions.  \hyperref[conj]{Conjecture \ref*{conj}} refines this conjecture as well.

One can generalize \cite[Theorem 3.3]{CSF} of Stanley and \cite[Theorem 5.3]{CQSF} of Shareshian and Wachs to show that the sums of certain coefficients in the $e$-basis expansion are positive.

\begin{prop}
Let $\overrightarrow{G}$ be a proper circular arc digraph on n vertices with G as its underlying undirected graph.  Suppose we have the expansion $X_{\overrightarrow{G}}({\bf x},t) = \displaystyle \sum_{\lambda \vdash n} c_{\lambda}(t)e_{\lambda}.$  Then $$\displaystyle \sum_{\substack{\lambda \vdash n \\ l(\lambda) = k}} c_{\lambda}(t) = \displaystyle \sum_{\bar{a} \in AO_k(G)} t^{\asc_{\overrightarrow{G}}(\bar{a})}$$ where $AO_k(G)$ is the set of acyclic orientations of $G$ with $k$ sinks and $\asc_{\overrightarrow{G}}(\bar{a})$ is the number of edges of $G$ for which $\bar{a}$ and $\overrightarrow{G}$ have the same orientation.
\end{prop}

The simplest proper circular arc digraphs that are not also unit interval digraphs are the cycles on n vertices, $\overrightarrow{C_n} := \overrightarrow{G}_{n,2}^*.$  In \cite[Proposition 5.4]{CSF}, Stanley provides a formula for the $e$-basis expansion of the chromatic symmetric functions of undirected cycles that show that they are $e$-positive.  We refine his formula for the chromatic quasisymmetric functions of directed cycles.

\begin{thm} \label{cycle expansion}
For directed cycles, $\overrightarrow{C_n},$ 
$$\displaystyle \sum_{n \geq 2} X_{\overrightarrow{C_n}}({\bf x},t)z^n = \frac{t \displaystyle \sum_{k \geq 2}k[k-1]_t e_k z^k}{1-t \displaystyle \sum_{k \geq 2} [k-1]_t e_k z^k}$$
where $[k]_t = 1 + t + \cdots + t^{k-1}.$ Consequently, $X_{\overrightarrow{C_n}}({\bf x},t)$ is palindromic, $e$-positive, and $e$-unimodal.
\end{thm}   

We prove\footnote{An alternative proof of this result was subsequently obtained in \cite{AP}.} this by extending the technique of Stanley for the $t=1$ case, which uses the transfer matrix method \cite{EC1}, and by using a result on permutation statistics by Mantaci and Rakotondrajao \cite{MR}. 

The following two propositions\footnote{These two propositions were obtained independently in \cite{AP}.} give combinatorial interpretations for the coefficients in the $e$-basis expansion of the chromatic quasisymmetric functions of $\overrightarrow{P_n},$ the path on n vertices oriented in one direction, and of $\overrightarrow{C_n},$ the cycle on n vertices oriented cyclically. The first proposition uses our generating function for the $e$-basis expansion of $X_{\overrightarrow{C_n}}({\bf x},t)$ seen in \hyperref[cycle expansion]{Theorem \ref*{cycle expansion}}. The second proposition uses a similar generating function by Shareshian and Wachs \cite{CQSF} for the $e$-basis expansion of $X_{\overrightarrow{P_n}}({\bf x},t).$  Though the proofs of these propositions rely on formulas that already give us $e$-positivity, perhaps there is a generalization of these interpretations that suggests a possible method for proving $e$-positivity for larger classes of graphs. 

\begin{prop}
Let $\overrightarrow{C_n}$ be the cycle on n vertices oriented cyclically with $C_n$ as its underlying undirected graph and let $X_{\overrightarrow{C_n}}({\bf x},t) = \displaystyle \sum_{\lambda \vdash n} c_{\lambda}(t) e_{\lambda}.$  Then $$c_{\lambda}(t) = \displaystyle \sum_{\bar{a} \in AO_{\lambda}(C_n)} t^{\asc_{\overrightarrow{C_n}}(\bar{a})}$$ where $AO_{\lambda}$ is the set of all acyclic orientations of $C_n$ such that the number of vertices between consecutive sinks of $\bar{a}$ is $\lambda_1 - 1,$ $\lambda_2 - 1,$...,$\lambda_k - 1$ in any order and $\asc_{\overrightarrow{C_n}}(\bar{a})$ is the number of edges of $\overrightarrow{C_n}$ for which the orientations of $\overrightarrow{C_n}$ and $\bar{a}$ agree.
\end{prop}

\begin{example}
Suppose we have the following acyclic orientation of $C_9,$ the underlying undirected graph of $\overrightarrow{C_9}.$  For convenience, we label the vertices with the elements of $[9]$ such that the edges of the original $\overrightarrow{C_9}$ are oriented from smaller label to larger label, except for the edge between 1 and 9, i.e. $1 \rightarrow 2 \rightarrow \cdots \rightarrow 9 \rightarrow 1.$

\begin{center}
\includegraphics[scale=0.4]{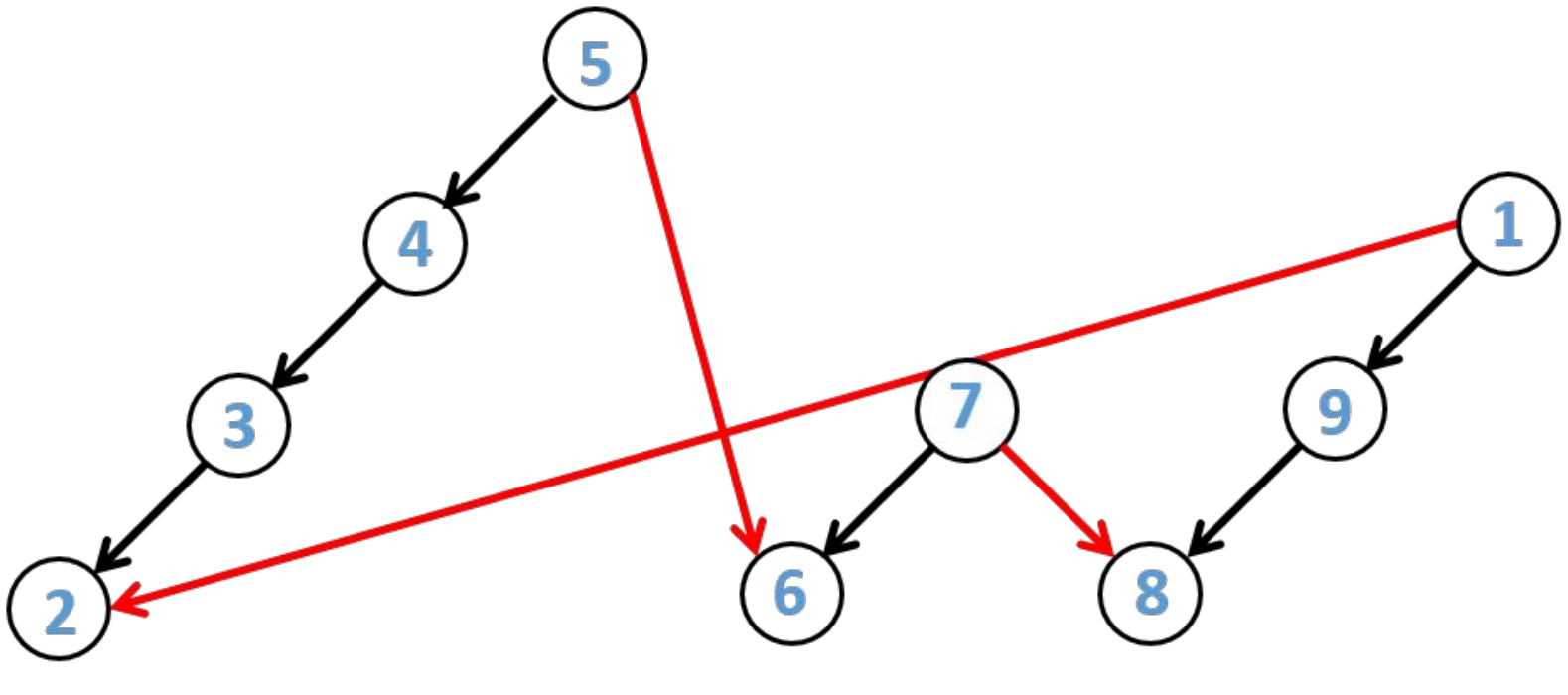}
\end{center}

This corresponds to a $t^3 e_{432}$ term, because there are 3 vertices between the sinks 2 and 6, 1 vertex between the sinks 6 and 8, and 2 vertices between the sinks 8 and 2.  The red edges are the 3 ascents of this orientation. 
\end{example}

\begin{prop}
Let $\overrightarrow{P_n}$ be the path of length n with edges oriented in one direction with $P_n$ as its underlying undirected graph and let $X_{\overrightarrow{P_n}}({\bf x},t) = \displaystyle \sum_{\lambda \vdash n} c_{\lambda}(t) e_{\lambda}.$  Then $$c_{\lambda}(t) = \displaystyle \sum_{\bar{a} \in AO_{\lambda}(P_n)} t^{\asc_{\overrightarrow{P_n}}(\bar{a})}$$ where $AO_{\lambda}(P_n)$ is the set of all acyclic orientations of $P_n$ such that the number of vertices between consecutive sinks of $\bar{a}$ (including the total number of vertices before the first sink and after the last sink) is $\lambda_1 - 1,$ $\lambda_2 - 1,$...,$\lambda_k - 1$ in any order and $\asc_{\overrightarrow{P_n}}(\bar{a})$ is the number of edges of $\overrightarrow{P_n}$ for which the orientations of $\overrightarrow{P_n}$ and $\bar{a}$ agree.
\end{prop}

\begin{example}
Suppose we have the following acyclic orientation of $P_8,$ the underlying undirected graph of $\overrightarrow{P_8}.$  For convenience, we label the vertices with the elements of $[8]$ such that the edges of the original digraph $\overrightarrow{P_8}$ are oriented from smaller label to larger label, i.e. $1 \rightarrow	2 \rightarrow \cdots \rightarrow 8$. 

\begin{center}
\includegraphics[scale=0.4]{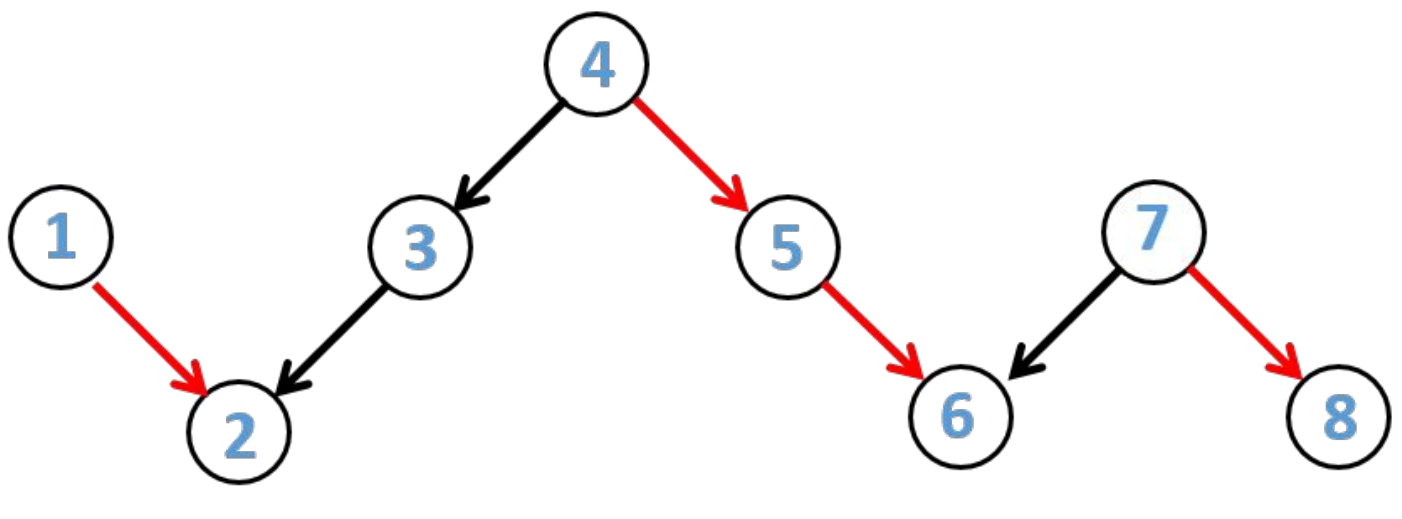}
\end{center}

This corresponds to a $t^4 e_{422}$ term, because there are 3 vertices between the sinks $2$ and $6$, 1 vertex between the sinks $6$ and $8$, 1 vertex before the 2 and none after the 8.  The red edges are the 4 ascents of this orientation. 
\end{example}



\section*{Acknowledgements}{I would like to thank Richard Stanley for his helpful suggestion that made this work possible.  I would also like to thank my advisor, Michelle Wachs, for all of her guidance and encouragement.}

\bibliography{FPSAC2017}
\bibliographystyle{alpha}

\end{document}